\documentclass[11pt]{article}
\usepackage{amsmath,amsfonts,amssymb}
\usepackage{amssymb}
\usepackage{epsfig}
\author{M.~Af\/fouf\\
School of Mathematics\\
Kean University\\
Union, NJ 07083 \\}
\date{ }
\title{A Remark on  Mathieu's Series }
\begin{document}
\maketitle

\begin{abstract}
We establish a new lower bound for Mathieu's series and present a new derivation of its expansions in terms of Riemann Zeta functions. \\
{\bf Key Words:}  Mathieu's series, Riemann Zeta functions\\
{\bf AMS(MOS) subject classifications.} 
Primary 26D15, Secondary 41A60.
\end{abstract}

In this work, we study Mathieu's series 
\begin{eqnarray}\label{1}  
F(h)= \sum_{n=1}^{\infty} \dfrac{n}{(n^2+h)^2} \medskip (h>0)
\end{eqnarray} 
which was derived in 1890 by Mathieu in his treatise on theory of elasticity and solid bodies \cite{4} and Mathieu conjectured the inequality \begin{eqnarray}\label{21}  
\sum_{n=1}^{\infty} \dfrac{2n}{(n^2+h)^2}< \dfrac{1}{2h} 
\end{eqnarray}
where $h $ is any real positive number.

In 1949, Schroder \cite{5} proved the the Mathieu's conjecture and another inequality:
 \begin{eqnarray}  \label{2} 
F(h)= \sum_{n=1}^{\infty} \dfrac{n}{(n^2+h)^2} < \dfrac{1}{2h}\medskip (h>0)
\end{eqnarray} 

 \begin{eqnarray}\label{3} 
F(h)< \dfrac{1}{(1+h)^2}+\dfrac{2}{(4+h)^2}+\dfrac{1}{2(4+h)} \medskip (0\leq h < 2)
\end{eqnarray} 


In 1960, Zmorovitch \cite{6} provided an elegant elementary proof of Schroder inequality based on the integral formula

\begin{equation}\label{22} 
\int_0^{\infty} e^{-ax} \cos bx dx = \dfrac{a}{a^2+b^2} ( a>0)
\end{equation}

Substitute $ a=n$ and $b=\sqrt{h}$, the last integral will have the form
\begin{equation}\label{4} 
\int_0^{\infty} e^{-nx} \cos \sqrt{h}x dx = \dfrac{n}{n^2+h} 
\end{equation}

Since the integral (\ref{4}) is absolutely convergent, we can take the derivative of both sides  with respect to the parameter $h$, to obtain

\begin{equation}\label{5} 
\int_0^{\infty} x  e^{-nx} \sin \sqrt{h}x dx = 2\sqrt{h} \dfrac{n}{(n^2+h)^2} 
\end{equation}

By taking the sum of both sides with respect to $n$ and using the geometric series formula, as follows:

\begin{equation}\label{6} 
\int_0^{\infty} x \sum_{n=1}^{\infty}  e^{-nx} \sin \sqrt{h}x dx = 2\sqrt{h} \sum_{n=1}^{\infty} \dfrac{n}{(n^2+h)^2} 
\end{equation}

\begin{equation}\label{7} 
\int_0^{\infty} \dfrac{x}{e^x-1}  \sin \sqrt{h}x dx = 2\sqrt{h} F(h) 
\end{equation}

We perform integration by parts three times to the integral on the left side of formula (\ref{7}), and we use the following expressions:

\begin{eqnarray}\label{8}
f(x)&=&(\dfrac{x}{e^x-1})'=\dfrac{1}{e^x-1}-\dfrac{xe^x}{(e^x-1)^2}\\ \nonumber
f'(x)&=&\dfrac{2xe^x}{(e^x-1)^3}-\dfrac{xe^x+2e^x}{(e^x-1)^2}\\ \nonumber
f''(x)&=&\dfrac{-6xe^{3x}}{(e^x-1)^4}-\dfrac{6xe^{2x}+6e^{2x}}{(e^x-1)^3}- \dfrac{xe^{x}+3e^{x}}{(e^x-1)^2} \nonumber
\end{eqnarray}

We obtain the formula

\begin{equation}\label{9} 
F(h)= \dfrac{1}{2h}+\dfrac{1}{2h^2}\int_0^{\infty} f''(x)(1-\cos \sqrt{h}x ) dx 
\end{equation}

Note that the function $f''(x)$ in the integrand of (9) can be written in the form 

\begin{equation}
f''(x)=\dfrac{(3-x)e^{2x}-4xe^{x}-x-3}{(e^x-1)^4}
\end{equation}

We can expand the numerator into the Taylor series such that

\begin{equation}
(3-x)e^{2x}-4xe^{x}-x-3=-\dfrac{x^5}{10}[1+x+\dfrac{23}{126}x^2+\dfrac{1}{14}x^3+\cdots]
\end{equation}
This implies that $f''(x) <0 $ and the Shroder inequality is directly obtained 

\begin{equation}
F(h) < \dfrac{1}{2h}
\end{equation}

Furthermore, By using L'hopital's  rule, we can find the limit

$$ \lim_{x\to 0} f'(x)=\dfrac{1}{6}$$
which implies that 
\begin{equation}
F(h) > \dfrac{1}{2h}-\dfrac{1}{6h^2}
\end{equation}
That is 

\begin{equation}
\dfrac{1}{2h}-\dfrac{1}{6h^2}< F(h)< \dfrac{1}{2h} 
\end{equation}
for all $h>0$.

Using the above discussed integral transformations, we establish the following proposition on expanding Mathieu's series in terms of Riemann Zeta functions.

{\bf Lemma:} The Mathieu's series can be expanded in terms of Riemann zeta functions
\begin{equation}
\sum_{n=1}^{\infty} \dfrac{n}{(n^2+h)^2}=
\sum_{n=1}^{\infty} \dfrac{(-1)^{n+1}}{(2n-1)!}\zeta(2n)h^{n-1}  
\end{equation}
for small parameter $h>0$.

{\bf Proof: }
We expand the sine function in  formula (7) in Taylor series for small parameter $h$,  as follows:
\begin{eqnarray}
2\sqrt{h} F(h)&=& \sqrt{h}
\int_0^{\infty} \dfrac{x}{e^x-1}[x-\dfrac{hx^3}{3!}+\dfrac{h^2x^5}{5!}-\dfrac{h^3x^7}{7!}-\dfrac{h^4x^9}{9!}+\cdots ]dx \\ \nonumber
2F(h)&=& \int_0^{\infty} \dfrac{x^2}{e^x-1}dx-
\dfrac{h}{3!}\int_0^{\infty} \dfrac{x^4}{e^x-1}dx+\dfrac{h^2}{5!}\int_0^{\infty} \dfrac{x^6}{e^x-1}dx \\ \nonumber
&-&\dfrac{h^3}{7!}\int_0^{\infty} \dfrac{x^8}{e^x-1}dx+\cdots \nonumber
\end{eqnarray}
These integrals are the Riemann zeta function so the new expression is

\begin{equation}
2F(h)=\zeta(2)-\dfrac{h}{3!}\zeta(4)+\dfrac{h^2}{5!}\zeta(6)-\dfrac{h^3}{7!}\varsigma(8)+\cdots
\end{equation}

From the last expression, we derive the new expansion of the Mathieu series in terms of Riemann zeta functions 

\begin{equation}
\sum_{n=1}^{\infty} \dfrac{n}{(n^2+h)^2}=
\sum_{n=1}^{\infty} \dfrac{(-1)^{n+1}}{(2n-1)!}\zeta(2n)h^{n-1}  
\end{equation}
or in an equivalent form

\begin{equation}
\sum_{n=1}^{\infty} \dfrac{n}{(n^2+h)^2}=
\dfrac{\pi^2}{6}-\dfrac{\pi^4}{90\cdot 3!}h+\dfrac{\pi^6}{945\cdot 5!}h^2-\cdots
\end{equation}

Clearly, for small $h$ the last formula implies that

$$ F(h)> \dfrac{\pi^2}{6}$$ which is better than available estimates. The proof is complete  $\blacksquare$.

Many questions remain to explain about the Mathieu series, here two of them: 


\begin{enumerate}
\item Investigate the generalized Mathieu series 
\begin{equation}
\sum_{n=1}^{\infty} \dfrac{n}{(n^2+h)^{\mu}}
\end{equation}
for $\mu>1$.

 \item Investigate the alternating Mathieu series 
\begin{equation}
S(h)=\sum_{n=1}^{\infty}(-1)^{n-1} \dfrac{n}{(n^2+h)^{2}}
\end{equation}
\end{enumerate}

Numerical investigation of the alternating Mathieu series, point to the fact $S(h)$ is a decreasing function with a range  $[0, 0.9015]$. However, this needs to be established analytically.




\end{document}